\documentclass{amsart}

\usepackage{amsmath,amssymb,amscd,amsthm}

\usepackage[latin1]{inputenc}
\usepackage{latexsym}
\usepackage{graphics}
\usepackage{ae}

\begin{document}

\def\HS{{\widetilde{\mathrm{HS}}}}

\title{Erratum to ``Polyhedral hyperbolic metrics on surfaces'' }

\author{Fran\c{c}ois Fillastre}

\date{September 11, 2008.}
\maketitle

In the last section of \cite{CompHyp} it is proved that the map $\mathcal{I}$ is a finite-sheeted covering map between 
$\mathcal{P}$ and $\mathcal{M}$. As $\mathcal{M}$ is  simply connected it is deduced that $\mathcal{I}$ is a homeomorphism. The fact that $\mathcal{P}$ is connected is missing. Here we provide a proof which is a simple adaptation of the argument for the case of the sphere  \cite[9.1.2]{AlexandrovBook}. The Fuchsian case with finite vertices was already done in \cite{Fillastre1}.

For the  parabolic case it is convenient to introduce a projective model of the hyperbolic-de Sitter space, which is a direct extension of a model of the hyperbolic space known as paraboloid model \cite[2.3.13]{Tlivre} or as parabolic model  \cite[10.25]{BriHa}. 
Let $\ell$ be a light-like vector of $\mathbb{R}^4_1$, the Minkowski space of dimension $4$. The paraboloid model can be obtained by a central projection in $\mathbb{R}^4_1$ of $\HS^3_{\ell}$ onto a light-like hyperplane parallel to $\ell^{\bot}$. In this model $\HS^3_{\ell}$ is homeomorphic to $\mathbb{R}^3$, geodesics are sent to straight lines, the boundary $\ell^{\bot}$ of $\mathrm{dS}^3_{\ell}$ is sent to infinity and we choose coordinates such that $\ell$ is sent to $(0,0,\infty)$. We denote by $D$ the part  of the boundary at infinity of $\mathbb{H}^3$  sent to the paraboloid  of equation $2z=-x^2-y^2$. The hyperbolic space is sent below $D$ and $\mathrm{dS}^3_{\ell}$ is sent above $D$. Horospheres centered at $\ell$ are sent to paraboloids $2z=-x^2-y^2+t$, where $t$ is a real number. The orthogonal projection onto the horospheres centered at $\ell$ along the lines starting from $\ell$ corresponds to the vertical projection  in this model.

For the Fuchsian case we consider the Klein projective model, the  plane fixed by the Fuchsian groups is the horizontal one and we take for $D$   the upper half part of the unit sphere. The orthogonal projection onto caps corresponds to the vertical projection  in this model.

 Let $(P,G)\in\mathcal{P}$. We denote by $\mathcal{P}^{(P,G)}$ the subset of  $\mathcal{P}$ made with elements $(Q,G)$  such that $Q$ has its vertices on the same vertical lines as the vertices of $P$. We also denote by $\overline{\mathcal{P}}^{(P,G)}$ the corresponding subset in the space of invariant convex HS-polyhedra (\emph{i.e.} we remove the condition that edges must meet hyperbolic space).

 The \emph{height} of a vertex $v$ is the signed vertical Euclidean distance between $D$ and $v$. For an ideal vertex the height is zero. It is positive for a strictly hyperideal vertex and negative for a finite vertex.  For example in the parabolic case  if $v$ has coordinates $(a,b,c)$, its height is $c+(a^2+b^2)/2$. 
Each element of $\overline{\mathcal{P}}^{(P,G)}$ is defined by the $(n+m+p)$ heights of vertices which are inside a fundamental domain for the action of $G$ (the polyhedral surface is then the convex hull of the orbits for $G$ of these $(n+m+p)$ vertices). In the following we identify $\overline{\mathcal{P}}^{(P,G)}$ with a subset of $\mathbb{R}^{n+m+p}$.

The condition for an element of $\mathbb{R}^{n+m+p}$ to correspond to an element of $\overline{\mathcal{P}}^{(P,G)}$ is that each vertex must lie outside the convex hull of the other vertices.
That means that if $v$ is one of the $(n+m+p)$ vertices and  $v_1,v_2,v_3$ are any other vertices  such that  $v$ is contained inside the vertical prism defined by $v_1,v_2,v_3$, then the plane spanned by $v_1,v_2,v_3$ must lie below $v$. This gives (an infinite number of) affine conditions on the heights of the $(n+m+p)$ vertices. 
We then get a convex subset of  $\mathbb{R}^{n+m+p}$ and if $(x_1,\ldots,x_{n+m+p})$ are the coordinates of the Euclidean space, we have to intersect this convex subset with $x_i<0$ if the $i$th coordinate corresponds to a finite vertex, with $x_i>0$ if it corresponds to a strictly hyperideal vertex and with $x_i=0$ if it is a finite vertex. Hence $\overline{\mathcal{P}}^{(P,G)}$ is a convex subset  of  $\mathbb{R}^{n+m}$.

 Let us call $(P^i,G)$ the ideal convex invariant polyhedron whose vertices are the vertical projection of the ones of $P$ onto $D$ ($(P^i,G)$ is the origin of $\mathbb{R}^{n+m}$ for the coordinates we introduced). It is lying on the boundary of $\overline{\mathcal{P}}^{(P,G)}$ as slightly pushing up and down suitable vertices of $(P^i,G)$   gives an element  of  $\overline{\mathcal{P}}^{(P,G)}$.

Hence in $\overline{\mathcal{P}}^{(P,G)}$ there is a segment $s$ from 
$(P,G)$ to $0$. Actually $s$ is lying in $\mathcal{P}^{(P,G)}$: along $s$ the heights of the strictly hyperideal vertices decrease, and if we start from a polyhedron with all edges meeting hyperbolic space, this property is preserved all along the deformation.

Let us denote by $\mathcal{P}^i$ the set of ideal convex invariant polyhedra with $(n+m+p)$ vertices in a fundamental domain. 
We described a continuous path in $\mathcal{P}$ from any element of $\mathcal{P}$ such that the other endpoint of its completion is in  $\mathcal{P}^i$. This last space  is path-connected as it is homeomorphic to the Teichm\"uller space of $(n+m+p)$ marked points on the compact surface $\overline{S}$.
Moreover any continuous path in $\mathcal{P}^i$ can be approximated by a path in $\mathcal{P}$.
Hence  $\mathcal{P}$ is path connected.

\bibliographystyle{alpha}

  \end{document}